\documentclass[a4paper,11pt]{article}
\usepackage{geometry}
\usepackage[dvipdfmx]{graphicx}
\usepackage{amsfonts}
\usepackage{amsmath}
\usepackage{url}
\begin{document}
\title{Performance evaluation of multiple precision matrix multiplications using parallelized Strassen and Winograd algorithms}
\author{Tomonori Kouya\thanks{Shizuoka Institute of Science and Technology}\\\url{http://na-inet.jp/}}
\date{2015-10-26}
%
%
\maketitle

\abstract{It is well known that Strassen and Winograd algorithms can reduce the computational costs associated with dense matrix multiplication. We have already shown that they are also very effective for software-based multiple precision floating-point arithmetic environments such as the MPFR/GMP library. In this paper, we show that we can obtain the same effectiveness for double-double (DD) and quadruple-double (QD) environments supported by the QD library, and that parallelization can increase the speed of these multiple precision matrix multiplications. Finally, we demonstrate that our implemented parallelized Strassen and Winograd algorithms can increase the speed of parallelized LU decomposition.}

\section{Introduction}
Multiple precision floating-point arithmetic has much heavier calculations for multiplication and division than for addition and subtraction. Therefore, Strassen\cite{strassen_original} and Winograd\cite{Coppersmith1990251}  algorithms, which can reduce the number of required multiplications in real matrix multiplication, are effective in the multiple precision floating-point environment. We have already implemented these two algorithms with MPFR (GNU MPFR) \cite{mpfr} and GMP (GNU MP) \cite{gmp}, and have shown their effectiveness for LU decomposition with matrix multiplication\cite{kouya_strassen2014}. In this paper, we show that the same effectiveness can be obtained in de facto standard quadruple precision (double-double, DD) and octuple precision (quadruple-double, QD) environments supported by Bailey's QD library \cite{qd}, and that parallelization using OpenMP can increase the speed of these algorithms. Naturally, these parallelized algorithms (DD, QD, and MPFR/GMP) can be effective for LU decomposition.

The QD library, developed by Bailey et al., implemented nearly quadruple precision (106 bits, DD) and octuple precision (212 bits, QD) floating-point numbers using two or four connected double precision floating-point numbers. The basic arithmetic is written in the C++ language; the library also provides ANSI C APIs. DD and QD floating point arithmetic operations are different from IEEE quadruple and octuple precision standards. Because the implementation scheme of the QD library is simple, many studies that use DD- and QD-type floating-point arithmetic have been performed in various areas of science and engineering.

On the other hand, MPFR is an integer-based arbitrary precision floating-point arithmetic library that provides IEEE754 floating-point standard compatible functions. MPFR is based on the multiple precision natural number (MPN) kernel supported by GMP; therefore we use the term ``MPFR/GMP" (MPFR over GMP) in this paper. Naturally, MPFR/GMP's performance depends on the MPN kernel of GMP, which has been well-tuned on various CPU architectures over the past 20 years. There are other multiple precision floating-point libraries that depend on GMP's MPN kernel, but MPFR/GMP is the best and oldest one.

We have been developing a multiple precision numerical computation library, BNCpack\cite{bnc}, based on MPFR/GMP; however, it does not use the QD library. In the rest of this paper, DD and QD precision matrix multiplications are implemented in the form of linear computation frameworks for BNCpack with MPFR/GMP.

\section{Strassen and Winograd algorithms and their parallelizations}

We consider a real matrix multiplication of any size defined as $C := AB = [c_{ij}]$ $\in\mathbb{R}^{m\times n}$, where $A =[a_{ij}]$ $\in \mathbb{R}^{m\times l}$ and $B = [b_{ij}]$$\in \mathbb{R}^{l\times n}$. Each element $c_{ij}$ of $C$ is defined as:
\begin{equation}
	c_{ij} := \sum^l_{k = 1} a_{ik} b_{kj}. \label{eqn:matrix_mul_simple}
\end{equation}
We call this simple algorithm ``Simple"; it uses the matrix multiplication presented in formula (\ref{eqn:matrix_mul_simple}) above.

We always use the ``block algorithm" (``Block"), which divides $A$ and $B$, to increase the hit ratio of the cache memory in the processor. In this paper, we divide $A$ and $B$ into small $M L$ pieces of $A_{ik}$, and $L N$ pieces of $B_{kj}$, respectively. Therefore, we can obtain blocked $C_{ij}$ using the following matrix multiplication:
\[
	C_{ij} := \sum^L_{k = 1} A_{ik} B_{kj}. 
\]

The complexity of block algorithm is the same as that of the simple algorithm. On the other hand, Strassen and Winograd algorithms can reduce complexity by employing recursive self calls.

\begin{figure}[htbp]
\begin{center}
\includegraphics[width=.7\textwidth]{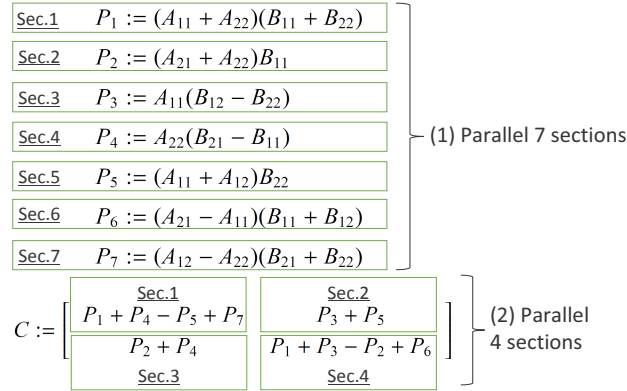}
\caption{Parallelized Strassen algorithm}\label{fig:parallel_strassen}
\end{center}
\end{figure}

\begin{figure}
\begin{center}
\includegraphics[width=.7\textwidth]{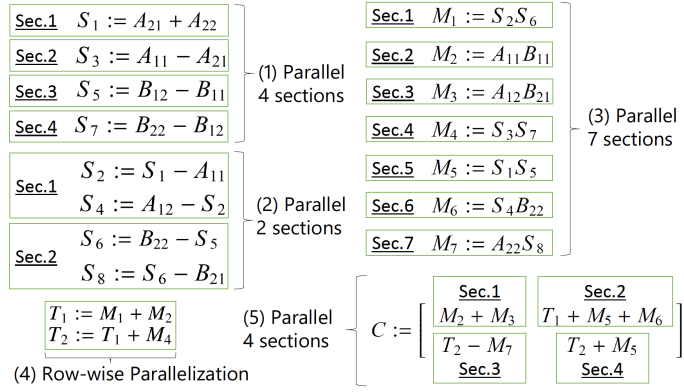}
\caption{Parallelized Winograd algorithm}\label{fig:parallel_winograd}
\end{center}
\end{figure}

As mentioned above, we have already implemented Strassen and Winograd algorithms using MPFR/GMP; the results are published as the BNCmatmul library\cite{bncmatmul}. These two algorithms are successfully able to shorten computational time, relative to block algorithms\cite{kouya_strassen2014}. To improve speed, first, we replaced all matches of ``\verb|omp parallel for|" in every loop; however, we cannot obtain more effective results. Second, we changed the parallelized Strassen and Winograd algorithms as shown in \figurename\ref{fig:parallel_strassen} and \figurename\ref{fig:parallel_winograd}. Independent and parallelizable parts of these algorithms are divided using ``\verb|omp| \verb|section|", and therefore, they can simultaneously execute their sections for each thread. The parallelized Winograd algorithm is more complex than the Strassen algorithm. Part (3) of the self recursive calls is divided into 7 threads. Likewise, part (1) of the parallelized Strassen algorithm is divided into 7 threads.

Because of these changes of parallelization, we can increase the speed of the Block, Strassen, and Winograd algorithms for matrix multiplications in not only MPFR/GMP but also DD and QD precision floating-point environments.

%
\section{Benchmark tests of parallelized matrix multiplication}
%
Our implemented C++ and C programs include parallelized ``Simple", ``Block(block size)", ``Strassen($n_{min}$)", and ``Winograd($n_{min}$)" algorithms with QD and MPFR/GMP libraries, where ``block size" refers to the minimal size of $A_{ik}$, and $B_{kj}$ and $n_{min}$ indicates the minimal size for which the recursive algorithm stops the self call. We have executed benchmark tests on the following computational environment.

\begin{description}
	\item[H/W] Intel Xeon E5-2620 v2 (2.10 GHz), 32 GB RAM
	\item[S/W] CentOS 6.5 x86\_64, Intel C/C++ 13.1.3, MPFR 3.1.2\cite{mpfr} / GMP 6.0.0a\cite{gmp} $+$ BNCpack 0.8\cite{bnc}, qd 2.3.15\cite{qd}
\end{description}
In the above environment, DD basic arithmetic, provided by the original QD library, is roughly 3 to 5 times faster than MPFR/GMP (106-bit mantissa) and QD is slightly slower than MPFR/GMP (212 bits). The QD library has the potential to increase in speed when using various techniques such as applying SIMD commands or FMA with assembler coding, which is applied to the well-tuned MPN library of GMP. We use original C++ classes such as \verb|dd_real| and \verb|qd_real| provided by the QD library.

The real square matrices $C:= AB$, $A$, $B$$\in \mathbb{R}^{n\times n}$ are as follows:
\[ A = \left[\sqrt{5} \left(i + j - 1\right)\right]^{n}_{i,j=1},\ B = \left[\sqrt{3}\left(n - i\right)\right]^{n}_{i,j=1} \]

\begin{table}[htb]
\begin{center}
\caption{Computational time of DD (nit: seconds)}\label{table:dd}
\begin{tabular}{|c|c|c|c|c|c|c|c|}\hline
1PE	& \multicolumn{4}{|c|}{DD (C++)} \\ \hline
$n$	& Simple	& Block(32)	& Strassen(32) & Winograd(32) \\ \hline
1023	& 32.3		& 20.8	& 11.7  & 11.7 \\
1024	& 49.6		& 20.3	& 11.7  & 11.6 \\
1025	& 32.6		& 22.3	& 11.7  & 11.7 \\ \hline
8PEs	& \multicolumn{4}{|c|}{DD (C++)} \\ \hline
$n$	& Simple	& Block(32)	& Strassen(32) & Winograd(32) \\ \hline
1023	& 68.3		& 3.2	& 3.2   & 3.2 \\
1024	& 69.7		& 3.2	& 3.1   & 3.1 \\
1025	& 68.3		& 4.1	& 3.2   & 3.2 \\ \hline
\end{tabular}
\end{center}
\end{table}

\begin{figure}[htb]
\begin{center}
\includegraphics[width=.7\textwidth]{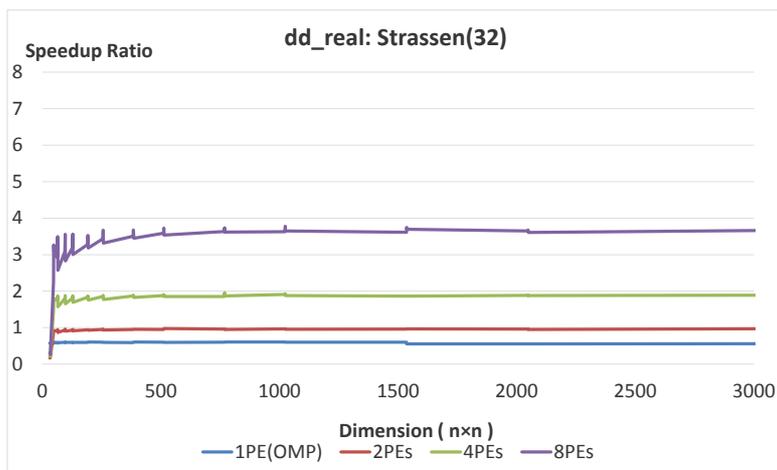}
\caption{Speed increase ratio of parallelized DD Strassen algorithm}\label{fig:strassen_parallel_dd_ratio}
\end{center}
\end{figure}

\tablename\ \ref{table:dd} shows the computational time of the DD square matrix multiplications. Parallelization of the Simple algorithm is not completely effective, but parallelized Block algorithms with DD and QD arithmetic see speed improvements (proportional to the number of threads). As shown in \figurename\ref{fig:strassen_parallel_dd_ratio}, parallelization of the Strassen algorithm with 8 threads can be up to 4 times faster than serial computation. The parallelized Winograd algorithm is at the same level of speed improvement as the Strassen algorithm.

\begin{figure}[htb]
\begin{center}
\includegraphics[width=.7\textwidth]{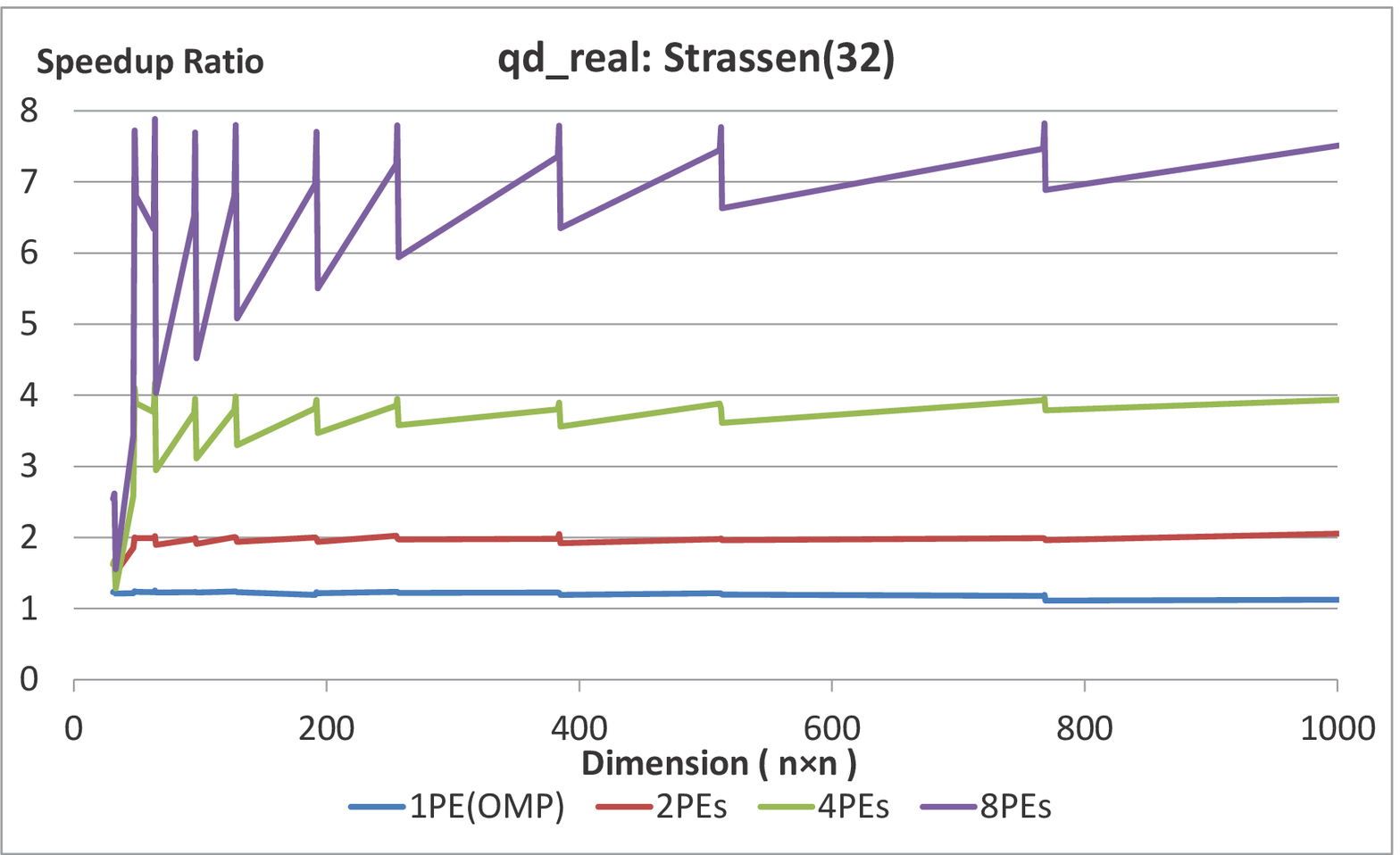}
\includegraphics[width=.7\textwidth]{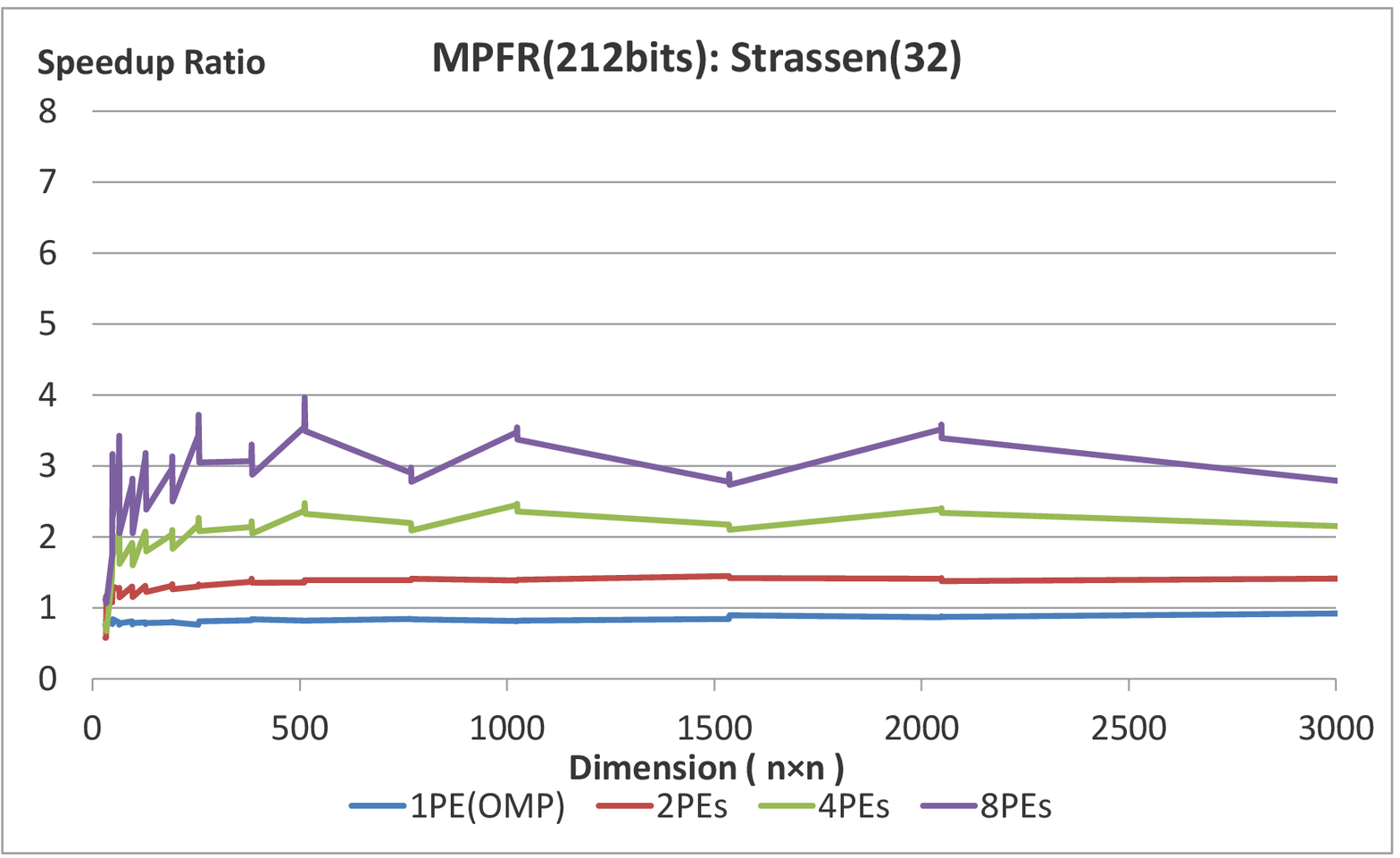}
\caption{Speed increase ratio of parallelized QD (upper) and MPFR/GMP (212 bits)(lower) Strassen algorithm}\label{fig:strassen_parallel_qd_ratio}
\end{center}
\end{figure}

\tablename\ \ref{table:qd_vs_mpfr} shows the computational time of QD and MPFR/GMP (212 bits) square matrix multiplications. Parallelization of the QD matrix multiplication with Block and Strassen algorithms can increase the speed, as shown in \figurename\ref{fig:strassen_parallel_qd_ratio}. Serial QD matrix multiplication is slower than MPFR/GMP (212 bits), as shown in \tablename\ \ref{table:qd_vs_mpfr}; however, parallel MPFR/GMP (212 bits) matrix multiplications are slightly slower than those of QD, because of the lower speed increase ratio of MPFR matrix multiplication.

\begin{table}[htb]
\begin{center}
\caption{Computational times of QD and MPFR/GMP (212 bits) matrix multiplication (unit: seconds)}\label{table:qd_vs_mpfr}
\begin{tabular}{|c|c|c|c|c|c|c|c|}\hline
1PE	& \multicolumn{3}{|c|}{QD (C++)} & \multicolumn{3}{c|}{MPFR (212bits)} \\ \hline
$n$	& B(32)	& S(32)	& W(32)	& B(32)	& S(32)	& W(32)	\\ \hline
1023	& 249.0	& 134.5	& 134.7	& 160.5	& 76.0	& 75.7 	\\
1024	& 247.6	& 134.3	& 134.5	& 163.2	& 75.1	& 75.2 	\\
1025	& 272.4	& 135.0	& 134.9	& 161.1	& 76.7	& 75.9 	\\ \hline
8PEs	& \multicolumn{3}{|c|}{QD (C++)} & \multicolumn{3}{c|}{MPFR (212bits)} \\ \hline
$n$	& B(32)	& S(32)	& W(32)	& B(32)	& S(32)	& W(32)\\ \hline
1023	& 32.5	& 17.8	& 17.8	& 23.5	& 21.9	& 21.9 \\
1024	& 32.6	& 17.2	& 17.2	& 23.5	& 21.2	& 20.9 \\
1025	& 42.8	& 18.9	& 18.9	& 28.0	& 22.7	& 23.2 \\ \hline
\end{tabular}
\end{center}
\end{table}

As all of the results of our benchmark tests show, DD, QD, and MPFR/GMP matrix multiplications exhibit increased speed when using parallelization with OpenMP.

\section{Application to parallelized LU decomposition}
%
Finally, we consider the application of parallelized algorithms so as to increase the speed of parallel LU decomposition when using DD, QD, and MPFR/GMP arithmetic. It is well known that matrix multiplication can be applied to LU decomposition \cite{golub4thed}. In this implementation, none of the LU decomposition involves pivoting.

We consider the linear equation (\ref{eqn:linear_eq}) with $A \in\mathbb{R}^{n\times n}$, $\mathbf{b}\in\mathbb{R}^n$:
\begin{equation}
	A\mathbf{x} = \mathbf{b}. \label{eqn:linear_eq}
\end{equation}
We use direct methods for LU decomposition of the coefficient matrix by setting the block size to $K = \alpha n_{min}$. Then, LU decomposition with matrix multiplication (the underlined term) is as follows: 
\begin{enumerate}
	\item Divide $A$ into $A_{11}\in\mathbb{R}^{K\times K}$, $A_{12}\in\mathbb{R}^{K\times (n - K)}$, $A_{21}\in\mathbb{R}^{(n - K)\times K}$, and $A_{22}\in\mathbb{R}^{(n-K)\times (n-K)}$.
	\item Decompose $A_{11}$ into $L_{11} U_{11} (= A_{11})$, and then, transform $A_{12}$ into $U_{12}$ and $A_{21}$ into $L_{21}$.
	\item $A^{(1)}_{22} := A_{22} - \underline{L_{21} U_{12}}$
\end{enumerate}
After substituting $A := A^{(1)}_{22}$, repeat the above algorithm until $n - K \geq 0$.


We employ a random matrix as an example of a well-conditioned matrix and a Lotkin matrix as an example of an ill-conditioned matrix.
\begin{description}
	\item[Random Matrix] $a_{ij}$ is a random number in $[-1, 1]$.
	\item[Lotkin Matrix] $a_{ij} = \begin{cases}
		1 & (i = 1) \\
		1 / (i + j - 1) & (i \geq 2) \\
	\end{cases}$
\end{description}
The true solution is $\mathbf{x}$ $ = [0\ 1\ ...\ n-1]^T$; we set $\mathbf{b} := A\mathbf{x}$. The condition numbers $\|A\|_1\|A^{-1}\|_1$ of the random matrix and the Lotkin matrix for $n = 1024$ are $4.4\times 10^{6}$ and $4.3\times 10^{1576}$, respectively.  For the Lotkin matrix, we must use more than 5260 bits for $n = 1024$.

\begin{figure}[htb]
\begin{center}
\includegraphics[width=.7\textwidth]{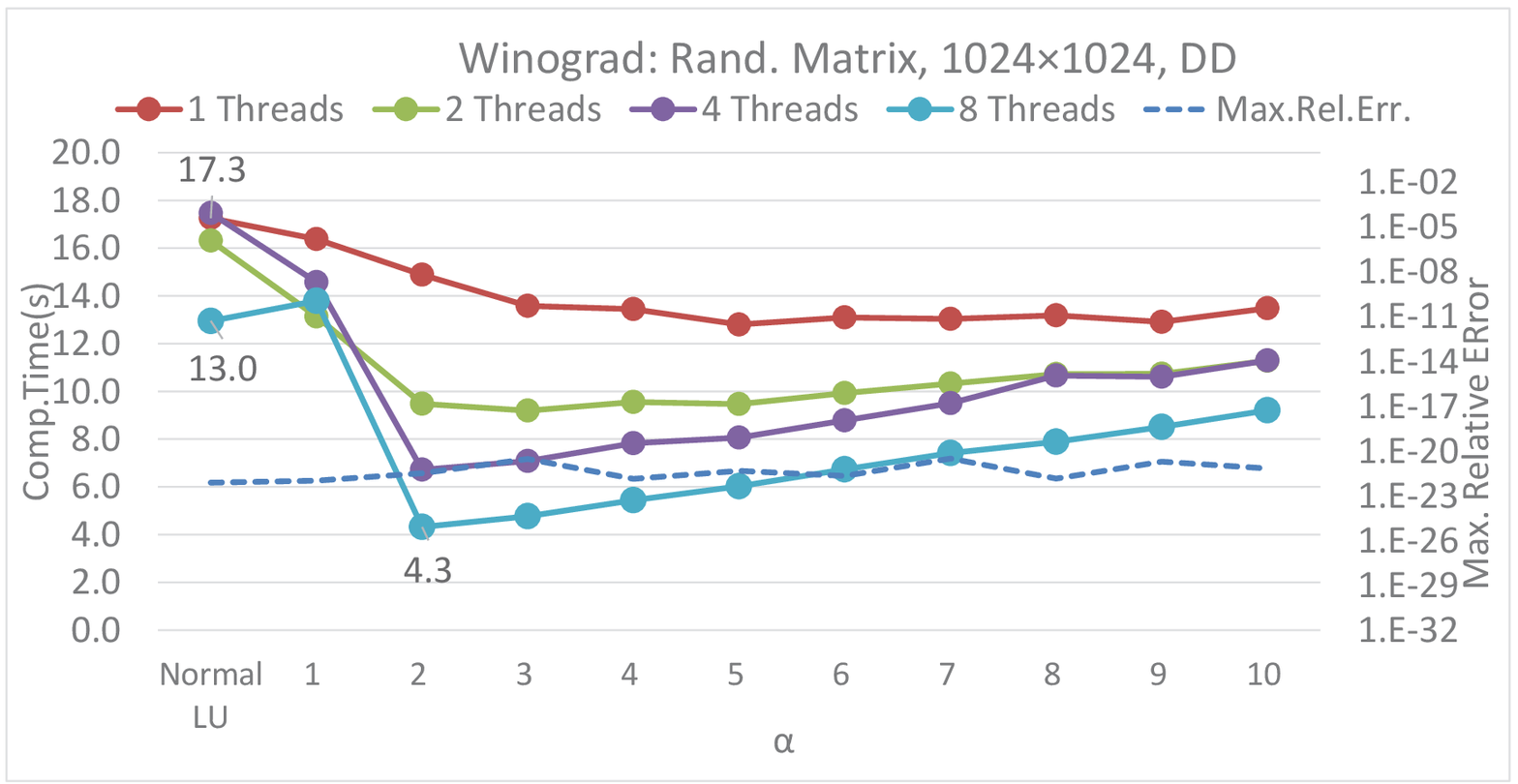}
\includegraphics[width=.7\textwidth]{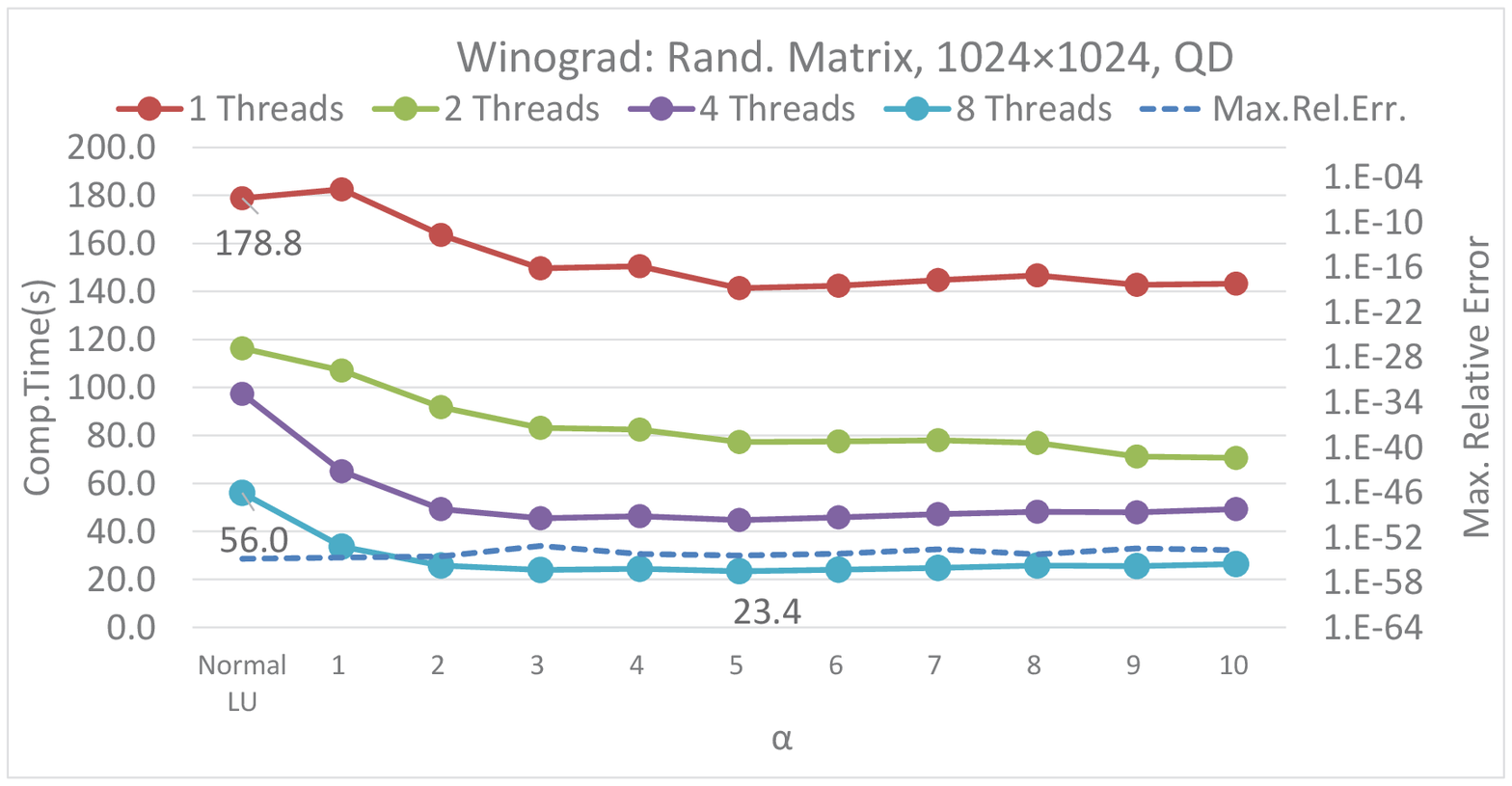}
\caption{Performance of DD (upper) and QD (lower) parallel LU decompositions and relative errors}\label{fig:parallel_lu_dd}
\end{center}
\end{figure}

\begin{figure}[htb]
\begin{center}
\includegraphics[width=.7\textwidth]{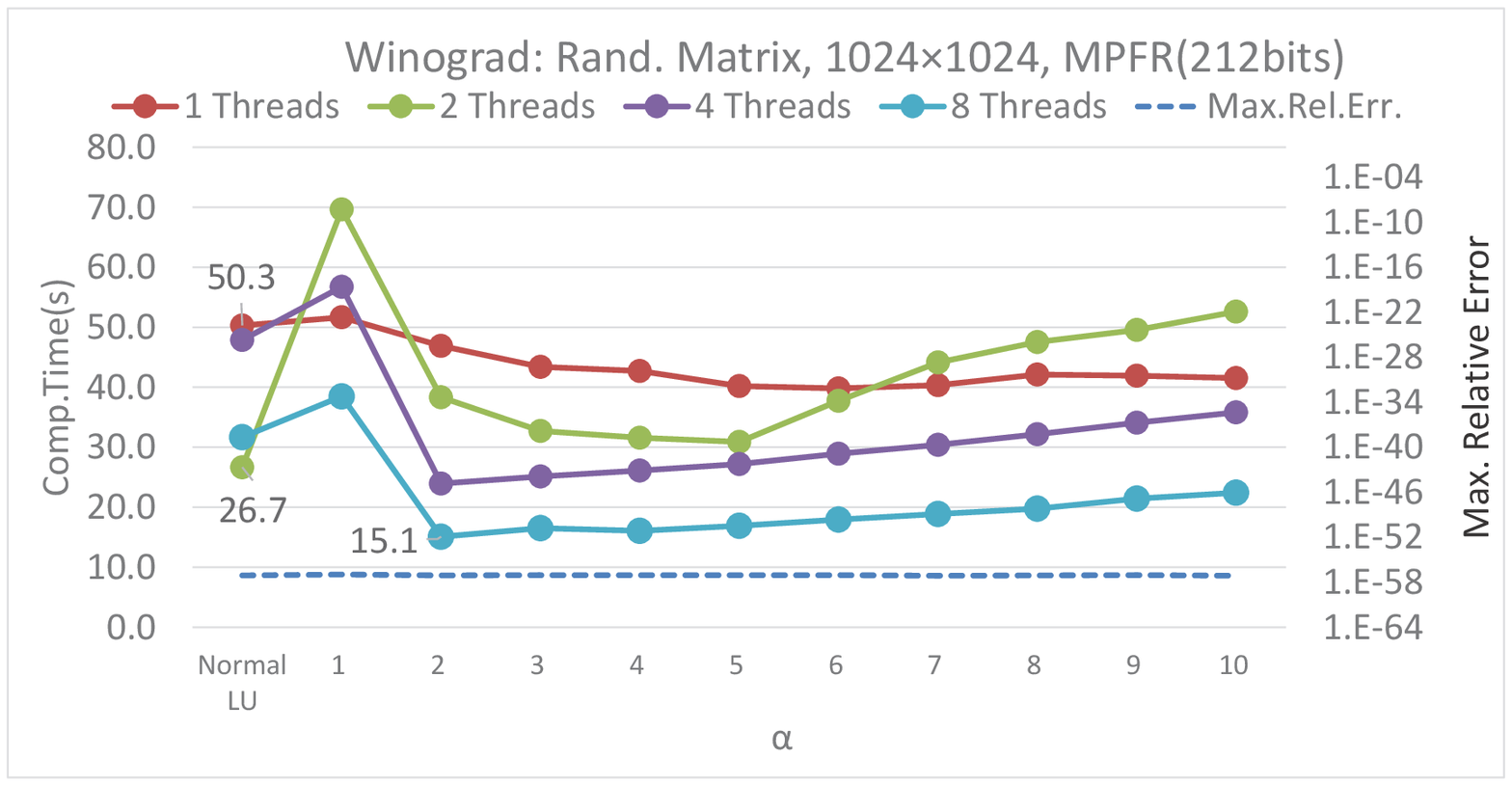}
\caption{Performance of MPFR/GMP (212 bits) parallel LU decomposition and relative error}\label{fig:parallel_lu_mpfr212}
\end{center}
\end{figure}

\begin{figure}[htb]
\begin{center}
\includegraphics[width=.7\textwidth]{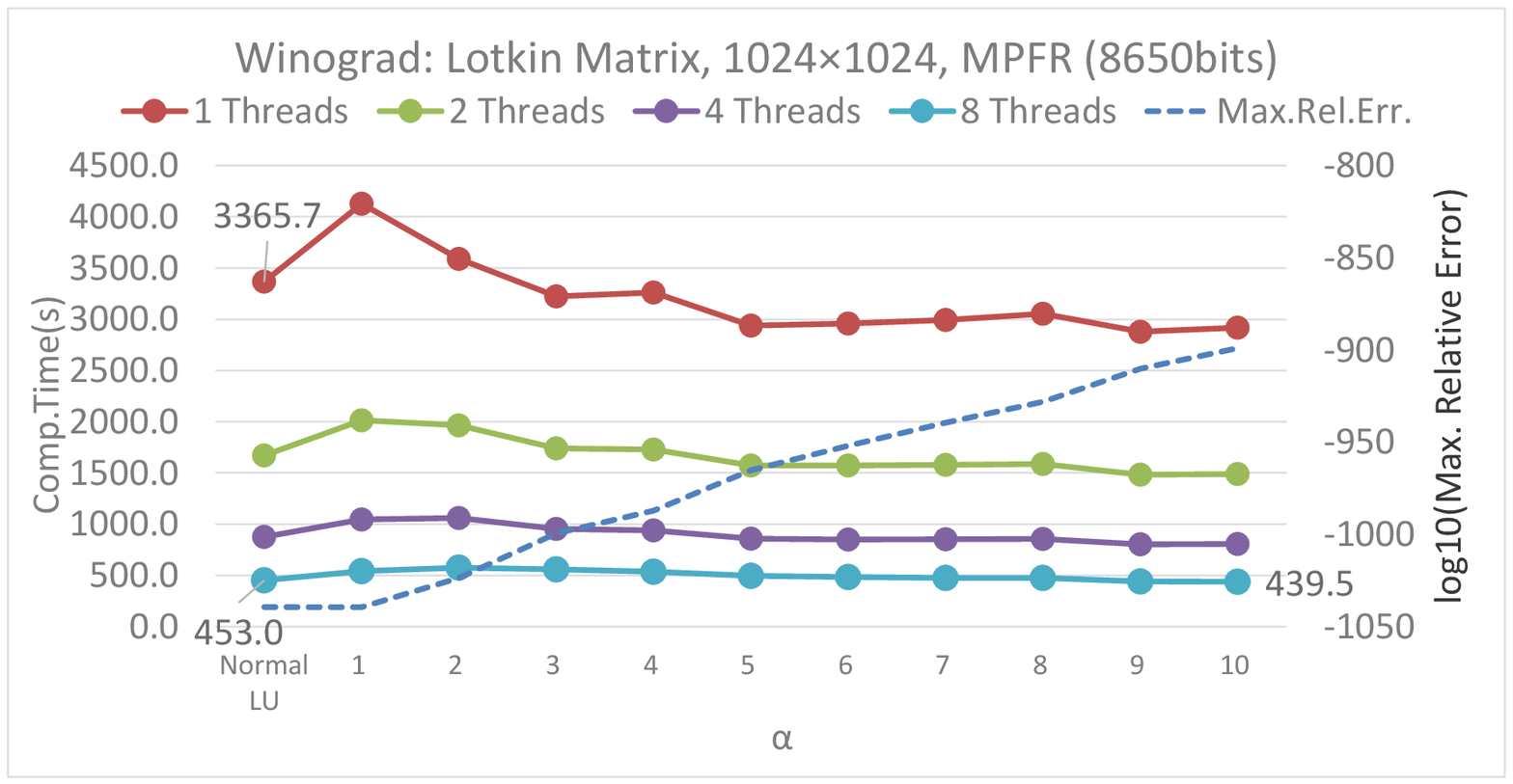}
\caption{Performance of MPFR/GMP (8650 bits) parallel LU decomposition and relative error}\label{fig:parallel_lu_mpfr8650}
\end{center}
\end{figure}

The $K$ sizes are set as $K = \alpha n_{min}$ ($\alpha = 1, 2, ..., 10$) and $n_{min} = 32$. Additionally, we investigated the computation time (in seconds) and the maximum relative error of the numerical solutions $\mathbf{x}$ at each value of $\alpha$. The random matrix is used with DD (\figurename\ref{fig:parallel_lu_dd} (upper)), QD (\figurename\ref{fig:parallel_lu_dd} (lower)), and MPFR/GMP (212 bits, \figurename\ref{fig:parallel_lu_mpfr212}), and the Lotkin matrix is used with MPFR/GMP (8650 bits, \figurename\ref{fig:parallel_lu_mpfr8650}). Because we know that the Winograd algorithm is slightly faster than the Strassen algorithm, as shown in the previous serial computation with MPFR/GMP \cite{kouya_strassen2014}, we compare the computational time using parallelized LU decomposition with the Winograd algorithm and row-wise parallelized LU decomposition.

We cannot determine the increase in relative errors that is due to the application of the Winograd algorithm in our benchmark tests for a random matrix. In the computational time of the DD LU decompositions, row-wise parallelized LU decomposition with 8 threads takes 13.0 s; however, the application of parallelized Wingorad algorithms can reduce it to 4.3 s at $\alpha = 2$ with 8 threads, as shown in \figurename\ref{fig:parallel_lu_dd}. QD LU decomposition can be reduced to 56 s with 8 threads when using row-wise parallelization; however, this reduces to 23.4 s at $\alpha = 5$ with 8 threads, as shown in \figurename\ref{fig:parallel_lu_dd}. In MPFR/GMP (212 bits) computation, we reduce the computation time from 26.7 s with 2 threads to 15.1 s at $\alpha = 2$ with 8 threads. Therefore, MPFR/GMP (212 bits) LU decomposition is slightly faster than QD decomposition.

For an ill-conditioned Lotkin matrix, we have already shown that Strassen and Winograd algorithms increase relative errors in the process of LU decomposition to 100 decimal digits (at most). Using parallelization, we can reduce computational time to 439.5 s at $\alpha = 10$ with 8 threads, as shown in \figurename\ref{fig:parallel_lu_mpfr8650}.

All of the results of our benchmark tests show that the parallelized Strassen and Winograd algorithms can increase the speed of multiple precision LU decomposition.

\section{Conclusions and future work}
Parallelization of Strassen and Winograd algorithms shows the following results.
\begin{itemize}
	\item Parallelized Strassen and Winograd algorithms as shown in \figurename\ref{fig:parallel_strassen} and \figurename\ref{fig:parallel_winograd} are effective on multicore CPUs; however, speed increase ratios for these algorithms are the same or lower than those of block algorithm.
	\item Serial MPFR/GMP (212 bits) matrix multiplication is faster than that of QD; however, parallelization reverses this tendency.
	\item Multiple precision LU decomposition can increase in speed when using Strassen and Winograd algorithms.
\end{itemize}
In our future work, we will increase the speed of multiple precision matrix multiplication on multicore CPUs and many-core environments such as GPUs and optimize both block size and times of self recursive calls in Strassen and Winograd algorithms.

%



\begin{thebibliography}{1}

\bibitem{qd}
D.H. Bailey, {QD}, 
\url{http://crd.lbl.gov/~dhbailey/mpdist/}.

\bibitem{Coppersmith1990251}
D. Coppersmith and S. Winograd,
Matrix multiplication via arithmetic progressions,
{\em Journal of Symbolic Computation}, {\bf 9}(1990), 251 -- 280.

\bibitem{golub4thed}
G. H. Golub and C. F. van Loan,
Matrix Computations (4th ed.),
Maryland,
Johns Hopkins University Press, 2013.

\bibitem{gmp}
T. Granlaud and {GMP} development team, 
The {GNU} {M}ultiple {P}recision arithmetic library,
\url{http://gmplib.org/}.

\bibitem{bnc}
T. Kouya, 
{BNC}pack, 
\url{http://na-inet.jp/na/bnc/}.

\bibitem{kouya_strassen2014}
T. Kouya, 
Accelerated multiple precision matrix multiplication using {S}trassen's algorithm and {W}inograd's variant,
{\em JSIAM Letters}, {\bf 6}(2014), 81 -- 84.

\bibitem{bncmatmul}
T. Kouya, {BNC}matmul, \url{http://na-inet.jp/na/bnc/bncmatmul-0.12.tar.bz2}

\bibitem{mpfr}
{MPFR} Project, 
The {MPFR} library, 
\url{http://www.mpfr.org/}.

\bibitem{strassen_original}
V. Strassen, 
Gaussian elimination is not optimal, 
{\em  Numerische Mathematik} {\bf 13}(1969), 354 -- 356.


\end{thebibliography}
\end{document}